\documentclass[11pt]{article}

\usepackage{amsmath}
\usepackage{amsthm}
\usepackage{amsfonts}
\usepackage{amssymb}
\usepackage{amscd}
\usepackage[dvips]{graphicx, color}
\usepackage[all]{xy}

\newtheorem{thm}{Theorem}[section]
\newtheorem{cor}[thm]{Corollary}
\newtheorem{defn}[thm]{Definition}
\newtheorem{lem}[thm]{Lemma}
\newtheorem{prop}[thm]{Proposition}
\newtheorem{rem}[thm]{Remark}

\newtheorem{con}[thm]{Conjecture}

\numberwithin{equation}{section}

\begin{document}

\title{\bf \LARGE On radiuses of convergence of $q$-metallic numbers and related $q$-rational numbers  }
\author{Xin Ren}

\date{}

\maketitle

\begin{abstract}
The $q$-rational numbers and the $q$-irrational numbers are introduced by S. Morier-Genoud and V. Ovsienko. In this paper, we focus on $q$-real quadratic  irrational numbers, especially $q$-metallic numbers and $q$-rational sequences which converge to $q$-metallic numbers, and consider the radiuses of convergence of them when we assume that $q$ is a complex number. We construct two sequences given by recurrence formula as a generalization of the $q$-deformation of Fibonacci numbers and Pell numbers which are introduced by S. Morier-Genoud and V. Ovsienko. We give an estimation of radiuses of convergence of them, and we solve on conjecture of the lower bound expected which is introduced by L. Leclere, S. Morier-Genoud, V. Ovsienko, and A. Veselov for the metallic numbers and its convergence rational sequence. In addition, we obtain a relationship between the radius of convergence of the $[n,n,\ldots,n,\dots]_q$ and $[n,n,\ldots,n]_q$ in the case of $n=3$ and $n=4$.
\end{abstract}


\section{Introduction}
\addcontentsline{toc}{section}{Introduction}

\noindent
It is well-known that the Euler $q$-integer is a polynomial with integer coefficients of $q$, which is a kind of quantization of integers. S. Morier-Genoud and V. Ovsienko expand integers to rational numbers and real numbers, and define $q$-rationals \cite{MO} and $q$-irrationals \cite{MO2} based on some combinatorial properties of rational numbers and number-theoretic properties of irrational numbers. The $q$-rationals and the $q$-irrationals are related to mathematical physics, combinatorics, number theory, quantum algebra, and knot theory, such as calculation recipes of the Jones polynomial and Alexander polynomial of rational knots (see \cite{KR},  \cite{KW} and \cite{NT}), quantum Teichm\"uller spaces \cite{VL}, and the Markov-Hurwitz approximation theory (see \cite{LMOV} and \cite{K}).

\medskip

We know that a positive rational number $\displaystyle \frac{r}{s}$ has regular and negative continued fraction expansions as $\displaystyle\frac{r}{s}=[a_1, \ldots , a_{2m}]=[[c_1, \ldots , c_k]]$. S. Morier-Genoud and V. Ovsienko \cite{MO} defined their $q$-deformations $\displaystyle [a_1, \ldots , a_{2m}]_q $ and $\displaystyle [[c_1, \ldots , c_k]]_q$, respectively. These $q$-deformations are equal, and thus one can define $\displaystyle \left[ \frac{r}{s} \right]_q:=
[a_1, \ldots , a_{2m}]_q = [[c_1, \ldots , c_k]]_q $ as $q$-rationals. On the other hand, the $q$-rationals can also be considered as a formal power series in $q$. Furthermore, given an arbitrary rational sequence $(x_k)_{k\geq1}$ that converges to an irrational number $x>1$,  one can consider the $q$-rational sequence $([x_k]_q)_{k\geq1}$. S. Morier-Genoud and V. Ovsienko proved that the Taylor series of $[x_k]_q$ stabilizes, as $k$ grows. This series does not depend on the choice of the sequence $[x_k]_q$. Hence, for an irrational $x>1$, the $q$-deformation of $x$ is defined as the following power series in $q$ : 
\[
[x]_q=\sum^{\infty}_{s=0}\varkappa_sq^s
\]
where $\displaystyle\varkappa_s$ $(s\in\mathbb{N})$ are integers \cite{MO2}.

\medskip

On the study of $q$-deformation of the real quadratic irrational numbers, L. Leclere and S. Morier-Genoud \cite{LM} proved some properties of $q$-deformation of them corresponding to classical real quadratic irrational numbers. On the other hand, since the $q$-real numbers are defined as power series in $q$, under the assumption that $q$ is a complex number, L. Leclere, S. Morier-Genoud, V. Ovsienko, and A. Veselov \cite{LMOV} study its radiuses of convergence and give the following conjecture which can be viewed as a $q$-deformation of Hurwitz's Irrational Number Theorem.
\begin{con}[\kern-0.35em{\cite{LMOV}}]\label{con-0.1}
{\normalfont
For any real number $x\geq1$, the radius of convergence of the series $[x]_q$ is greater than or equal to $\displaystyle R_{(1)}=\frac{3-\sqrt5}{2}$.
}
\end{con}

\indent
For convenience, we use the notation $\displaystyle R_{(1)}$ instead of $R_{\ast}$ in \cite{LMOV}.  It is shown that $R_{(1)}$ is the radius of convergence of the $q$-golden number $\displaystyle \left[\frac{1+\sqrt{5}}{2}\right]_q$. In the case where $x=[a_1,1,a_3,1,a_5,1,\ldots]$ ($a_1\geq1$, $a_{2m-1}\geq2$) and $x=[n,n,\ldots,n,\ldots]$ ($n=1,2$), Conjecture~\ref{con-0.1} has been proved in \cite{LMOV}. In this paper we solve Conjecture \ref{con-0.1} for the metallic numbers and its convergence rational sequence. More precisely, the following two theorems are proved:

\medskip

\begin{thm}\label{thm-0.2}
{\normalfont
For the metallic number $x=[n,n,\ldots,n,\ldots]$ with $n\geq3$, 
the radius of convergence of the series $[x]_q$ is greater than or equal to $\displaystyle R_{(1)}$.
}
\end{thm}

\medskip

\begin{thm}\label{thm-0.3.1}
{\normalfont
If $x=[n,n,\ldots,n]$, then the radius of convergence of the series $[x]_q$ is greater than 
or equal to $\displaystyle R_{(1)}$.
}
\end{thm}

\medskip

It is worth noting that the ideas for the proof of Theorems \ref{thm-0.2} and \ref{thm-0.3.1} are different. We convert the problem of finding the radius of convergence of the $q$-metallic number to consider the zeros of its discriminant, then the palindrome of the discriminant of the $q$-metallic number is the key to proving Theorem \ref{thm-0.2}. For the case of $[n,\ldots,n]_q$, its radius of convergence can also be converted to the zeros of the polynomial by considering the denominator part. On the other hand, S. Morier-Genoud and V. Ovsienko define $q$-deformation of Fibonacci numbers and Pell numbers (see \cite{MO}). Inspired by their study, we construct two sequences $\{\mathcal{M}_k(n)\}_{k=0,1,2,\ldots}$ and $\{\widetilde{\mathcal{M}}_k(n)\}_{k=0,1,2,\ldots}$ as a generalization of the $q$-deformation of Fibonacci numbers and Pell numbers. Different from the case of the $q$-metallic number, we mainly use the Rouch\'e theorem to complete the proof of Theorem \ref{thm-0.3.1} by considering the sequence $\{\mathcal{M}_k(n)\}_{k=0,1,2,\ldots}$.

\medskip

In addition, for some relationship between the radius of convergence of the $q$-metallic number $[n,n,\ldots,n,\ldots]_q$ and its truncated $q$-rational number $[n,n,\ldots,n]_q$, the case of $n=1,2$ has been examined in \cite{LMOV}. This means that for each $q$-rational in the $q$-rational sequence that converges to the $q$-silver number $[2,2,\ldots,2,\ldots]_q$, its radius of convergence is not only greater than or equal to $R_{(1)}$ but also greater than or equal to the radius of convergence of the $q$-silver number.  Based on this result, we observe the case of $n=3,4$, and give the following conclusion. 

\medskip

\begin{thm}\label{thm-0.3}
{\normalfont
For $n=3,4$, the radius of convergence of the series $[n,n,\ldots,n]_q$ is greater than the radius of convergence of the series $[n,n,\ldots,n,\ldots]_q$.
}
\end{thm}

\medskip

The idea of the proof is almost the same as the one of Theorem \ref{thm-0.3.1}. We still need to use the Rouch\'e theorem. However, for the general case ($n\geq5$), there are some limitations to this approach to find conditions that satisfy the Rouch\'e theorem. If this problem is solved, we can determine the relationship between the radiuses of convergence of  the $q$-metallic number $[n,n,\ldots,n,\ldots]_q$ and $[n,n,\ldots,n]_q$ in the general case.

\medskip

This paper is organized into the following sections.\\
\indent 
 In $\S 2$, we give definitions and some properties about $q$-rationals and $q$-irrationals introduced in \cite{MO} and \cite{MO2}, and we consider the $q$-deformed metallic number and its discriminant. 
  
\indent 
  In $\S 3$, we consider the radius of convergence of the $q$-metallic number. Firstly, we know that finding the radius of convergence can be transformed into finding the zeros of the discriminant. In particular, we deal with the discriminant appropriately to obtain a palindromic polynomial $P_n(q)$. We use its related properties to prove that $P_n(q)$ does not have a zero in the disk with radius $R_{(1)}$ and center at the origin. From this fact we prove Theorem \ref{thm-0.2}.
  
 \indent 
 In $\S 4$, we generalize the $q$-deformation of Fibonacci numbers and Pell numbers and define the recurrence formula of the sequence $\{\mathcal{M}_k(n)\}_{k=0,1,2,\ldots}$. We also give an important  property about them and prove Theorems \ref{thm-0.3.1}, \ref{thm-0.3} by this property. In this section, the Rouch\'e theorem plays an important role when we prove Theorems \ref{thm-0.3.1}, \ref{thm-0.3}, and related lemmas, we almost need to use it.

\section{The definition and some properties of $q$-real numbers}

\noindent
At the beginning, we briefly describe the definitions of the $q$-rational number and the $q$-irrational number. For details, it can be referred to \cite{MO} and \cite{MO2}. Then we consider the $q$-metallic number and its discriminant. The discriminant of the $q$-metallic number has palindromic properties.

\subsection{$q$-continued fractions and $q$-rational numbers}
\noindent
A positive rational number $\displaystyle\frac{r}{s}>1$ with coprime integers $r$ and $s\neq0$ has regular and negative continued fraction expansions as follows:\\
\[ \displaystyle
\frac{r}{s}=a_1+\frac{1}{a_2+\cfrac{1}{\ddots+\cfrac{1}{a_{n}}}}
=c_1-\frac{1}{c_2-\cfrac{1}{\ddots-\cfrac{1}{c_k}}}
\]
\\
with $a_i\geq1$ and $c_j\geq2$. These expansions are denoted by $[a_1,a_2,\ldots,a_{n}]$ and $[[c_1,c_2,\ldots,c_k]]$, respectively.
Since $[a_1,a_2,\ldots,a_{n}]=[a_1,a_2,\ldots,a_{n}-1,1]$, we always assume that $n=2m$ is even. In this way, we have the unique regular and negative continued fraction expansion of $\displaystyle\frac{r}{s}$.\\

\noindent 

Let $q$ be a formal parameter. For a non-negative integer $n$ we denote the Euler $q$-integer by 
\begin{equation}\label{eqn-1.1}
\displaystyle
[n]_q:=\left\{    
\begin{array}{
    cl}
0 & \text{for} \ n=0, \\
\displaystyle\frac{1-q^n}{1-q} =1 + q+q^{2} + \cdots + q^{n-1}  & \text{for} \ n\geq1. \\
\end{array} \right.
\end{equation}\
We note that $\displaystyle q^n[n]_{q^{-1}}=\frac{q^n(1-q^{-n})}{1-q^{-1}}=q[n]_q$.

\noindent 

The $q$-deformations of regular and negative continued fraction expansions are defined by Morier-Genoud and Ovsienko \cite{MO} as follows:
 
\begin{equation}\label{eqn-1.2}
[a_1, a_2, \ldots , a_{2m}]_q :=
[a_1]_q  
+\cfrac{ q^{a_1}   }{ [a_2]_{ q^{-1} }  + \cfrac{ q^{-a_2}  }
{  [a_3]_q  + \cfrac{q^{a_3} }
{ [a_4]_{q^{-1} } + \cfrac{ q^{-a_4} }{  \cfrac{\ddots }
{ [a_{2m-1}]_q+ \cfrac{ q^{a_{2m-1}} }{ [a_{2m}]_{q^{-1}}  }}}}}} 
\end{equation}

\noindent 

\begin{equation}\label{eqn-1.3}
[[c_1, c_2, \ldots , c_{k}]]_q :=
[c_1]_q  
-\cfrac{ q^{c_1 -1}   }{ [c_2]_{ q }  - \cfrac{ q^{c_2 -1}  }
{  [c_3]_q  - \cfrac{q^{c_3-1 }}
{ [c_4]_{q } - \cfrac{ q^{c_4 -1} }{  \cfrac{\ddots }
{ [c_{k-1}]_q- \cfrac{ q^{c_{k-1}-1} }{ [c_{k}]_{q}  }}}}}} . 
\end{equation}

\begin{thm}[{Morier-Genoud and Ovsienko \cite{MO}}] \label{thm-1.1}{\normalfont
If a positive rational number $\displaystyle\frac{r}{s}$ is given in the form 
\noindent 
$\displaystyle\frac{r}{s}=[a_{1}, \ldots , a_{2m}]=[[c_{1}, \ldots , c_{k}]]$, then 

\medskip

\noindent 
\begin{equation}\label{eqn-1.4}
[a_{1}, \dots , a_{2m}]_{q}=[[c_{1}, \ldots , c_{k}]]_{q}.
\end{equation}
Hence, the $q$-rational number $\displaystyle\left[\frac{r}{s}\right]_q$ is defined as 
\[\displaystyle\left[\frac{r}{s}\right]_q:=[a_{1}, \ldots , a_{2m}]_{q}=[[c_{1}, \ldots , c_{k}]]_{q}.\]
}
\end{thm}

\medskip



\subsection{$q$-irrationals}
Let $ x >1 $ be an irrational number. Let $ (x_k) _ {k \geq1} $ be a rational number sequence  that converges to $ x $, and consider the sequence $ ([x_k] _q)_{ k \geq1} $ of $q$-deformation. For $ k \geq1 $, we express $ [x_k] _q $ as the formal power series :\\
\begin{equation}\label{eqn-1.5}
\displaystyle [x_k]_q=\sum_{s=0}^{\infty}\varkappa_{k,s}q^s.
\end{equation}
Then the $q$-deformed irrational number $x$ is defined as the following formal power series in $q$:
\begin{equation}\label{eqn-1.6}
\displaystyle
[x]_q=\sum_{s=0}^{\infty}\varkappa_{s}q^s \ \ \ \ (\varkappa_s=\lim_{k \to \infty}\varkappa_{k,s})
\end{equation} 
We can guarantee that the existence of the limit and its independence of the choice of the converging sequence $(x_k)_{k\geq1}$ by the following proposition.
\begin{prop}[{Morier-Genoud and Ovsienko \cite[Theorem 1]{MO2}}] {\normalfont
Given an irrational real number $x\geq1$, for every $k\geq0$, the coefficients $\varkappa_{k,s}$ of (\ref{eqn-1.5}) are stabilizing, as $k$ grows. Moreover, the limit $\varkappa_s$ of coefficients(\ref{eqn-1.6}) are integers, and independent of the choice of the converging sequence $(x_k)_{k\geq1}$.}
\end{prop}
By the above proposition, we may always construct $q$-irrational numbers using continued fractions as follows:
for $x = [a_1,a_2,a_3,\ldots]>1$, the sequence of rational numbers which converges to $x$ can be chosen as 
$
\displaystyle
x_k=[a_1,a_2,\ldots,a_k].
$

\subsection{$q$-metallic numbers and its discriminant}
\begin{defn}[{$q$-metallic numbers}] {\normalfont
Let $n$ be a natural number. The metallic number $x$ can be represented as follows :\[\displaystyle x=\displaystyle \frac{n + \sqrt{n ^ {2} + 4}} {2} =[n,n,\ldots,n,\ldots].\] 

\indent
We call that the $q$-deformation of $x$ (denoted by $[x]_q$) the ($n$th) $q$-metallic number. Particularly, when $n=1$, the $[x]_q$ is called the $q$-golden number (the $1$st $q$-metallic number).
}
\end{defn}

\medskip

By (\ref{eqn-1.2}), we can obtain the next formula.
\begin{equation}\label{eqn-1.4.1}
\displaystyle
\left[\frac{n+\sqrt{n^2+4}}{2}\right]_q =
[n]_q  
+\cfrac{ q^{n}   }{ [n]_{ q^{-1} }  + \cfrac{ q^{-n}  }
{  [n]_q  + \cfrac{q^{n} }
{ [n]_{q^{-1} } + \cfrac{ q^{-n} }{\ddots }
}}}
\end{equation}
The $q$-metallic number is a solution of the algebraic equation
\begin{equation}\label{eqn-1.4.2}
\displaystyle
q[x]_q^{2}-((q-1)(q^n+1)+q[n]_q)[x]_q-1=0, 
\end{equation}
which is a $q$-deformation of $x^2-nx-1=0$. This can be deduced from (\ref{eqn-1.4.1}).
\medskip

\begin{defn}{\normalfont
We denote by $\mathcal{D}([ \overline{n}]_q)$ the discriminant which is satisfied as the quadratic equation (\ref{eqn-1.4.2}) for $[x]_q$. More precisely, it is given by
\begin{equation}
\displaystyle
\mathcal{D}([ \overline{n}]_q)
=\left\{    
\begin{array}{
    cl}
1+2q-q^2+2q^3+{q}^{4} & \text{for} \ n=1, \\
1+2q+3q^2+3q^4+2q^5+{q}^{6}  & \text{for} \ n=2, \\
D(q) & \text{for} \ n\geq3, \\
\end{array} \right.
\end{equation}
where
\[
\displaystyle
D(q)=1+2q^2+\sum^{n-1}_{t=3}(t-1)q^t+(n+1)q^n+(n-4)q^{n+1}+(n+1)q^{n+2}+\sum^{2n-1}_{t=n+3}(2n-t+1)q^t+2q^{2n}+q^{2n+2}.
\]
Then $\mathcal{D}([\overline{n}]_q)$ is said to be the discriminant of the $q$-metallic number.
}
\end{defn}
\indent
We note that all of the coefficients of $\mathcal{D}([ \overline{n}]_q)$ are positive integers for $n\geq4$. The discriminant $\mathcal{D}([\overline{n}]_q)$ has the following property. This property is very helpful for us to prove Theorem \ref{thm-0.2}.

\medskip

\begin{thm} {\normalfont(the palindromic of discriminant)
\\
The discriminant of $q$-metallic numbers $\mathcal{D}([ \overline{n}]_q)$ are palindromic. 
}
\end{thm}

\noindent
Proof. It is obvious that $\mathcal{D}([ \overline{1}]_q)$ and $\mathcal{D}([ \overline{2}]_q)$ are palindromes. For $n\geq3$, one can check $q^{2n+2}\mathcal{D}([ \overline{n}]_{q^{-1}})=\mathcal{D}([ \overline{n}]_q)$ by direct calculation.
\qed

\medskip

\begin{rem}{\normalfont
For general $q$-real quadratic irrational numbers, the palindromic of discriminant has been proved by considering a $q$-deformation of the elements of the modular group $PSL(2,\mathbb{Z})$ (see \cite{LM}). 
}
\end{rem}

\section{The radius of convergence of $q$-metallic numbers}
Since $q$-irrational is defined as power series in $q$, we can consider its radiuses of convergence. In this section, we always assume that $q$ is a complex number.
\subsection{Some observations}
For $q$-metallic numbers \[\displaystyle [x]_q=\left[\displaystyle \frac{n + \sqrt{n ^ {2} + 4}} {2} \right]_q=[n,n,\ldots,n,\ldots]_q,\] with $n\geq1$,
we denote the radiuses of convergence of $[x]_q$ by $R_{(n)}$. In the cases of $n=1$ and $n=2$ the following has already been known.

\medskip

\begin{prop}[{\cite[Propositions 3.1 and 4.1]{LMOV}}]\label{prop-2.1.1} {\normalfont
\[
\displaystyle
R_{(1)}=\frac{3- \sqrt{5}}{2}=1-\frac{1}{\frac{1+ \sqrt{5}}{2}}\approx0.38197, \ R_{(2)}=\frac{1+\sqrt{2}- \sqrt{2\sqrt{2}-1}}{2}\approx0.53101.
\]
}
\end{prop}

\medskip

In complex analysis, it is well-known that the radius of convergence of a power series $[x]_q$ centered on a point $0$ is equal to the distance from $0$ to the nearest point so that $[x]_q$ cannot be defined in a way that makes it holomorphic. Since
 \[\displaystyle [x]_q=\frac{q[n]_q+(q^n+1)(q-1)+\sqrt{\mathcal{D}([ \overline{n}]_q)}}{2q},\] we can calculate $R_{(n)}$ by finding the roots of $\mathcal{D}([ \overline{n}]_q)$. 

\noindent

Since the discriminant $\mathcal{D}([ \overline{n}]_q)$ is a palindromic polynomial, we can find the roots of it by using the next lemma. 

\begin{lem}[\kern-0.35em{\cite[Corollary 2]{KM}}]\label{lem-2.1} {\normalfont Let $P(q) \in \mathbb{Q}[q]$ be a palindromic polynomial given by 
\[
P(q)=a_0+a_1q+a_2q^2+\cdots+a_{n-2}q^{n-2}+a_{n-1}q^{n-1}+a_nq^n
,\]
having even degree $n$. If
\[
\displaystyle
\max\{\left|a_k\right|: k\in {0,1,2,\ldots,\frac{n}{2}-1}\} \geq \left|a_{\frac{n}{2}}\right|,
\]  
then $P(q)$ has roots on the unit circle.}
\end{lem}

\medskip 

Since $\mathcal{D}([ \overline{n}]_q)$ satisfies Lemma \ref{lem-2.1}. it has roots on the unit circle. Indeed, for any positive integer $n$, we have
\[
\displaystyle
\mathcal{D}([ \overline{n}]_q)=(1-q+q^2)P_n(q),
\]
where 
\[
\displaystyle
P_n(q):=\left\{    
\begin{array}{cl}
1+3q+q^2 &  \text{for}  \   n=1, \\
1+q+4q^2+q^3+q^4  & \text{for}  \ n=2, \\
\displaystyle 1+\sum^{n-1}_{t=1}tq^t+(n+2)q^{n}+\sum^{2n-1}_{t=n+1}(2n-t)q^t+q^{2n} & \text{for}  \ t=n\geq3. \\
\end{array} \right.
\] 
Hence, we obtain $R_{(n)}\leq 1$, and $R_{(n)}=\min\{|q| : P_n(q)=0\}$.

\medskip 

\begin{lem}\label{lemma-2.4} {\normalfont
A polynomial $P(q)\in \mathbb{Q}[q]$ is palindromic if and only if the following two conditions are satisfied:
\begin{enumerate}
\item[(i)] $1$ is a root of even multiplicity (possibly zero),
\item[(ii)] if $q \neq \pm 1$ is a root, then $\displaystyle \frac{1}{q}$ is a root with the same multiplicity.
\end{enumerate}
}
\end{lem}
For a proof of Lemma \ref{lemma-2.4}, see \cite[Theorem 1.1.9]{L}.

\medskip

\subsection{Proof of Theorem \ref{thm-0.2}. }
We prove Theorem \ref{thm-0.2} by reductive absurdum. Let $n\geq3$, and assume that $q=re^{i\theta}$ ($0<r\leq R_{(1)}$, $0\leq\theta\leq2\pi$) is a root of $P_n(q)$. By Lemma~\ref{lemma-2.4}, since $P_n(q)$ is palindrome, it follows that $q^{-1}=r^{-1}e^{-i\theta}$ is also a root of $P_n(q)$. Hence, $P_n(q)$ has a root can be expressed by $q=rz$ ($\displaystyle r\geq \frac{1}{R_{(1)}}>\frac{5}{2}$, $|z|=1$). Then, $z$ is a zero of the polynomial

\begin{equation}\label{eqn-2.2.1}
\displaystyle
\sum^{2n}_{t=0}\mu_tz^t, \ \ 
\displaystyle
\mu_t=\left\{    
\begin{array}{
    cl}
1 &  \text{for}  \   t=0, \\
tr^t  & \text{for}  \ 0<t<n, \\
(n+2)r^n & \text{for}  \ t=n, \\
(2n-t)r^t & \text{for}  \  n<t<2n, \\
r^{2n} & \text{for}  \ t=2n. \\
\end{array} \right.
\end{equation}

\medskip

\begin{lem}\label{Lem2.2.1} 
{\normalfont For $n\geq3$, the sequence $\{\mu_t\}^{2n}_{t=0}$ with $\mu_t>0$ for $t=0,1,2,\ldots, 2n$, is strictly increasing and $\displaystyle \mu_{2n-1}<\frac{2}{5}\mu_{2n}$.}
\end{lem}

\noindent
Proof. Since $\displaystyle r>\frac{5}{2}$, then $\displaystyle \mu_{2n-1}=r^{2n-1}<\frac{2}{5}r^{2n}=\frac{2}{5}\mu_{2n}$. \qed

\medskip

We are now turning to the proof of Theorem \ref{thm-0.2}. Let us consider the polynomial \[\displaystyle F(z)=(1-z)\sum^{2n}_{t=0}\mu_tz^t=-\mu_{2n}z^{2n+1}+(\mu_{2n}-\mu_{2n-1})z^{2n}+\cdots+(\mu_1-\mu_0)z+\mu_0.\]
By Lemma \ref{Lem2.2.1} and $|z|=1$, one has
\[
\begin{split}
\displaystyle 
\left|F(z)\right|&=\left|-\mu_{2n}z^{2n+1}+\mu_{2n}z^{2n}-\frac{2}{5}\mu_{2n}z^{2n}+(\frac{2}{5}\mu_{2n})z^{2n}+\cdots+(\mu_1-\mu_0)z+\mu_0\right| \\
&\geq \left|\mu_{2n}\right|\left|z\right|^{2n}\left|z-\frac{3}{5}\right|-\left|z\right|^{2n}\left|(\frac{2}{5}\mu_{2n}-\mu_{2n-1})+\frac{(\mu_{2n-1}-\mu_{2n-2})}{z}+\cdots+\frac{(\mu_1-\mu_0)}{z^{2n-1}}+\frac{\mu_0}{z^{2n}}\right| \\
&\geq \left|\mu_{2n}\right|\left|z\right|^{2n}\left|z-\frac{3}{5}\right|-\left|z\right|^{2n}(\left|\frac{2}{5}\mu_{2n}-\mu_{2n-1}\right|+\frac{\left|\mu_{2n-1}-\mu_{2n-2}\right|}{\left|z\right|}+\cdots+\frac{\left|\mu_1-\mu_0\right|}{\left|z\right|^{2n-1}}+\frac{\left|\mu_0\right|}{\left|z\right|^{2n}}) \\
&= \mu_{2n}(\left|z-\frac{3}{5}\right|-\frac{2}{5}) = r^{2n}(\left|z-\frac{3}{5}\right|-\frac{2}{5})\geq0 .
\end{split}
\]
Note that the equality holds if and only if $z=1$, but $z=1$ is not a zero of (\ref{eqn-2.2.1}). It turns out that (\ref{eqn-2.2.1}) has no zero on the unit circle. This leads to a contradiction.\qed

\medskip

By the part (ii) of Lemma \ref{lemma-2.4} and the proof of Theorem \ref{thm-0.2}, it follows that $P_n(q)$ has no zeros inside $\displaystyle \{q\in\mathbb{C}:|q|>\frac{1}{R_{(1)}}=\frac{3+\sqrt{5}}{2}\}$. Then we have the following corollary. 

\begin{cor}\label{cor2.2.1}{\normalfont
All zeros of $P_n(q)$ lie in the annulus $\displaystyle \mathcal{A}_1:=\{q\in\mathbb{C}:\frac{3-\sqrt{5}}{2}\leq|q|\leq\frac{3+\sqrt{5}}{2}\}$. 
}
\end{cor}

\medskip

\begin{rem}\label{rem-2.2.1}{\normalfont
We calculate $R_{(n)}$ ($n=1,2,\ldots,48$) by using computer program (Maple) and correct to 6 decimal places (see, Table \ref{table:001}). From Table \ref{table:001}, we can know that $R_{(n)}$ does not monotonous increase such as $n=20,32,33,45,\ldots$. 
\begin{table}[h]
 \caption{the radiuses of convergence of $q$-metallic numbers}
 \label{table:001}
 \centering
  \begin{tabular}{clclclcl}
   \hline
   $n$ & $R_{(n)}$ & $n$ & $R_{(n)}$ & $n$ & $R_{(n)}$& $n$ & $R_{(n)}$\\
   \hline \hline
   1 & 0.38197 & 13 &  0.84033 & 25 &0.89505 & 37 &0.91965 \\
   2 & 0.53101 & 14 &  0.84047 & 26 &0.89483 & 38 &0.91970 \\
   3 & 0.59719 & 15 &  0.84280 & 27 &0.89505 & 39 &0.91952 \\
   4 & 0.65060 & 16 &  0.84767 & 28 &0.89665 & 40 &0.92025 \\
   5 & 0.69918 & 17 &  0.85496 & 29 &0.89967 & 41 &0.92193 \\
   6 & 0.74444 & 18 &  0.86423 & 30 &0.90395 & 42 &0.92449 \\
   7 & 0.76933 & 19 &  0.87477 & 31 &0.90917 & 43 &0.92775 \\
   8 & 0.77406 & 20 &  \textcolor{red}{0.87404}& 32 &\textcolor{red}{0.90916} & 44 &0.92784 \\
   9 & 0.78191 & 21 &  0.87485& 33 &\textcolor{red}{0.90910} & 45 &\textcolor{red}{0.92759} \\
   10 & 0.79338 & 22 & 0.87747& 34 &0.91016 & 46 &0.92811 \\
   11 & 0.80802 & 23 & 0.88192& 35 &0.91236 & 47 &0.92944 \\
   12 & 0.82492 & 24 & 0.88793& 36 &0.91560 & 48 &0.93153 \\
   \hline
  \end{tabular}
\end{table}
}
\end{rem}

\medskip 

\section{$q$-rational number sequence with convergence on $q$-metallic numbers}

\subsection{Two sequences $\{\mathcal{M}_k(n)\}_{k=0,1,2,\ldots}$ and $\{\widetilde{\mathcal{M}}_k(n)\}_{k=0,1,2,\ldots}$} 
Such as the golden number, that is the $1$st metallic number, can be represented by the limit of the quotient of the adjacent two terms of the Fibonacci sequence, the $n$th metallic number also has a sequence of numbers that is the limit of the quotient of the adjacent two terms. In this section, we first consider such a sequence and its $q$-deformation.

\indent
Fix a natural number $n$, and consider the sequence $\{A_k(n)\}_{k=0,1,2,\ldots}$ defined by $A_0(n)=0$, $A_1(n)=1$, and $A_{k+2}(n)=nA_{k+1}(n)+A_{k}(n)$. The sequence $\{A_k(n)\}_{k=0,1,2,\ldots}$ satisfies the following recurrence formula: 
\begin{equation}\label{eqn-3.1.1}
A_{k+4}(n)=(n^2+2)A_{k+2}(n)-A_k(n) \ (k=0,1,2,\ldots).
\end{equation}
We consider the rational sequence $\frac{A_2(n)}{A_1(n)},\frac{A_3(n)}{A_2(n)},{\frac{A_4(n)}{A_3(n)},\ldots}$, which converges the $n$th metallic number, and their $q$-deformations ${\left[\frac{A_2(n)}{A_1(n)}\right]_q, \left[\frac{A_3(n)}{A_2(n)}\right]_q, \left[\frac{A_4(n)}{A_3(n)}\right]_q,\ldots}$.

 It is well-known that $A_k(1)$ is the $k$th Fibonacci number, and $A_k(2)$ is the $k$th Pell number. Since the $q$-deformation of Fibonacci numbers and Pell numbers are defined in \cite{MO}. S. Morier-Genoud and V. Ovsienko \cite{MO} introduced two sequences which converge to $q$-golden number and $q$-silver number. These sequence can also be viewed as a $q$-deformation of Fibonacci numbers and Pell numbers, respectively. Inspired by their study, we introduce two sequences which are generalizations of them as follows.

\begin{prop}\label{prop3.1.1}{\normalfont Consider the following four recurrence formulas: 
\begin{equation}\label{eqn-3.1.2a}
\displaystyle
\begin{split}
&\mathcal{M}_{2l+1}(n)=q[n]_q\mathcal{M}_{2l}(n)+\mathcal{M}_{2l-1}(n),\\
&\mathcal{M}_{2l+2}(n)=[n]_q\mathcal{M}_{2l+1}(n)+q^{2n}\mathcal{M}_{2l}(n), \\
\end{split}
\end{equation}
\begin{equation}\label{eqn-3.1.2b}
\begin{split}
\displaystyle
&\widetilde{\mathcal{M}}_{2l+1}(n)=[n]_q\widetilde{\mathcal{M}}_{2l}(n)+q^{2n}\widetilde{\mathcal{M}}_{2l-1}(n),\\
&\widetilde{\mathcal{M}}_{2l+2}(n)=q[n]_q\widetilde{\mathcal{M}}_{2l+1}(n)+\widetilde{\mathcal{M}}_{2l}(n), \\
\end{split}
\end{equation}
with $l=1,2,\cdots$, and \[\mathcal{M}_{0}(n)=0,\ \mathcal{M}_{1}(n)=1, \ \mathcal{M}_{2}(n)=[n]_q, \ \mathcal{M}_{3}(n)=1+q([n]_q)^2,\] \[\widetilde{\mathcal{M}}_{0}(n)=0, \ \widetilde{\mathcal{M}}_{1}(n)=1, \ \widetilde{\mathcal{M}}_{2}(n)=[n]_q, \ \widetilde{\mathcal{M}}_{3}(n)=([n]_q)^2+q^{2n-1}.\]
Then, we have
\begin{equation}\label{eqn-3.1.3}
\displaystyle
\left[\frac{A_{k+1}(n)}{A_k(n)}\right]_q=\frac{\widetilde{\mathcal{M}}_{k+1}(n)}{\mathcal{M}_k(n)} \ (k=1,2,3,\ldots).
\end{equation}
}
\end{prop}

\noindent
Proof. We show this proposition by using the $q$-deformed continued fraction and weighted triangulations. For a detailed definition of weighted triangulation, it can be referred to \cite{MO} (see, Sections 2.1 and 2.3). 

\noindent

Then $[n,n,\ldots,n,\ldots]_q$ corresponds to the weighted triangulation depicted in Figure 1, where the weights of the unmarked edges are $1$ and the convergence of the continued fraction corresponds to the black vertices. By direct deduction, the numerator and denominator of $\displaystyle\left[\frac{A_{k+1}(n)}{A_k(n)}\right]_q$ are satisfy (\ref{eqn-3.1.2a}) and (\ref{eqn-3.1.2b}), respectively.
\begin{figure}[htbp]
\centering
\includegraphics[height=7.0cm, width=15.7cm]
{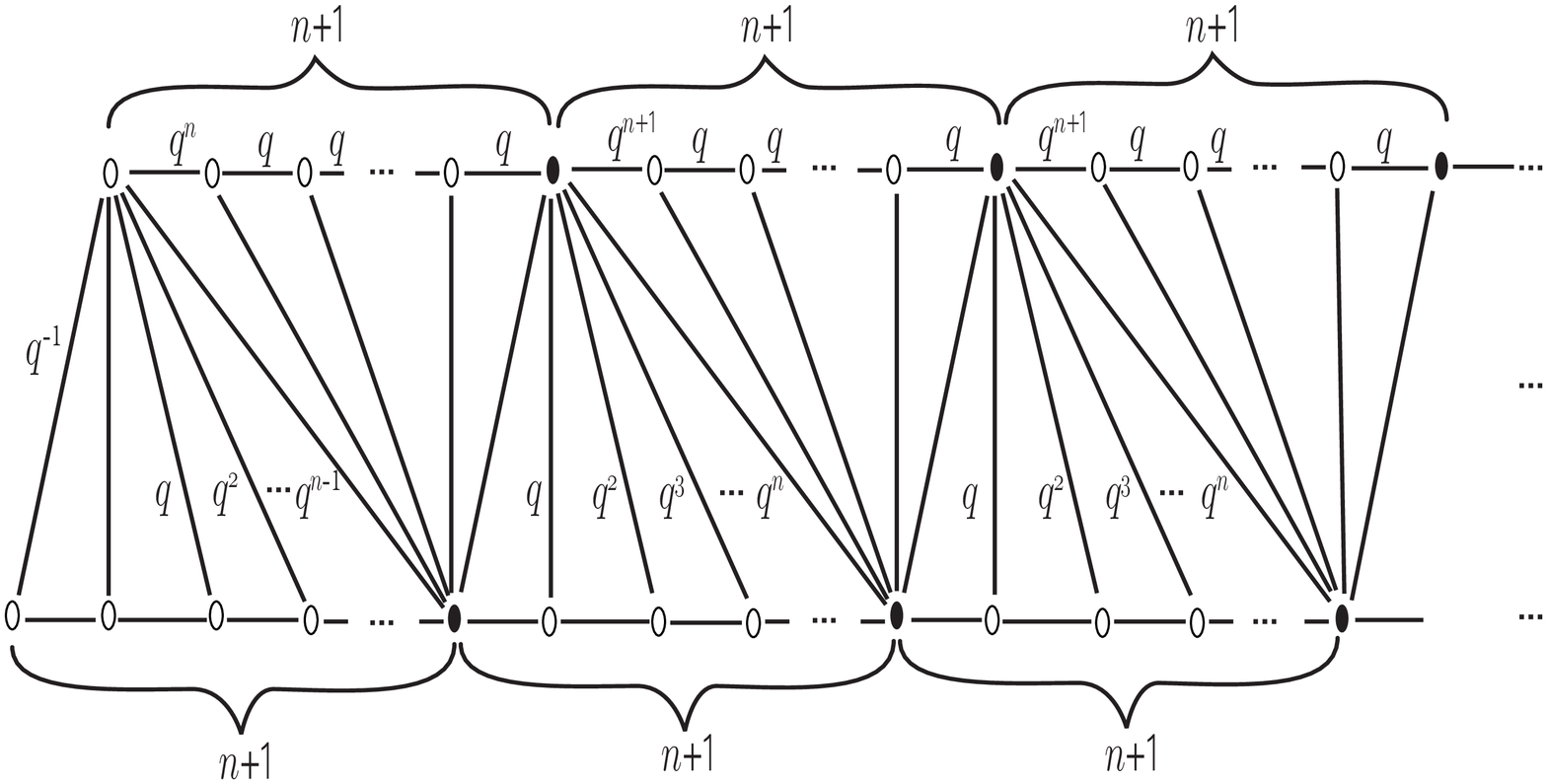}
\caption{
The weighted triangulation corresponding to $[n,n,\ldots,n,\ldots]_q$
}
\end{figure}
\qed

\medskip

Both $\{\mathcal{M}_k(n)\}_{k=0,1,2,\ldots}$ and $\{\widetilde{\mathcal{M}}_k(n)\}_{k=0,1,2,\ldots}$ can be viewed as the $q$-deformation of the sequence $\{A_k(n)\}_{k=0,1,2,\ldots}$. From (\ref{eqn-3.1.2a}) and Proposition \ref{prop3.1.1}, we can obtain the following proposition,  which can be viewed as a $q$-deformation of \eqref{eqn-3.1.1}.

\begin{prop}\label{prop3.1.2}{\normalfont The $\mathcal{M}_k(n)$ is determined by the recurrence formula
\begin{equation}\label{eqn-3.1.4}
\displaystyle
\mathcal{M}_{k+4}(n)=f(q,n)\mathcal{M}_{k+2}(n)-q^{2n}\mathcal{M}_{k}(n) \ (k=1,2,3,\ldots),
\end{equation}
where 
\[
\displaystyle
f(q,n):=1+q([n]_q)^2+q^{2n}=1+\sum^{n}_{t=1}tq^t+\sum^{2n-1}_{t=n+1}(2n-t)q^t+q^{2n},
\] 
and the initial data
\[\mathcal{M}_{1}(n)=1, \ \mathcal{M}_{2}(n)=[n]_q, \ \mathcal{M}_{3}(n)=1+q([n]_q)^2, \ \mathcal{M}_{4}(n)=(q^{2n}+1)[n]_q+q([n]_q)^{3}.\] 
}
\end{prop}
Note that $f(q,n)+2q^n=P_n(q)$, and so $f(q,n)$ is also a palindromic polynomial. 
\medskip


\subsection{Proof of Theorem \ref{thm-0.3.1}.}
For any $n$, the equation
\begin{equation}\label{eqn-3.2.1}
\displaystyle
\frac{\mathcal{M}_{k+2}(n)}{\mathcal{M}_{k}(n)}=f(q,n)-q^{2n}\frac{\mathcal{M}_{k-2}(n)}{\mathcal{M}_{k}(n)}
\end{equation} can be implied by \eqref{eqn-3.1.4}. We set $C_n=\{q\in\mathbb{C}:|q|=R_{(n)}\}$ and $D_n=\{q\in\mathbb{C}:|q|<R_{(n)}\}$. Before starting the proof, we will prepare 3 lemmas as follows. 


\begin{lem}\label{lem3.2.01}{\normalfont On the circle $C_{1}$, one has
\[
\displaystyle
\left|f(q,n)\right|_{C_1}>1-R_{(1)}.
\]
}
\end{lem}

\noindent
Proof. On the circle $C_1$, one has
\[
\displaystyle
\left|f(q,n)\right|_{C_1}\geq f(-R_{(1)},n)=1+\sum^{n}_{t=1}(-1)^{t}tR_{(1)}^{t}+\sum^{2n-1}_{t=n+1}(-1)^t(2n-t)R_{(1)}^{t}+R_{(1)}^{2n}.
\]
Furthermore, the $f(-R_{(1)},n)$ in the above inequality can also be written as
\[
\displaystyle 
f(-R_{(1)},n)=1+\sum_{\substack{t\in\left\{0,\ldots,n\right\}\\2\mathrel{|}t}}R_{(1)}^t(t-(t+1)R_{(1)})+\sum_{\substack{t\in\{n+1,\ldots,2n-2\}\\2\mathrel{|}t}}R_{(1)}^t((2n-t)-(2n-t-1)R_{(1)})+R_{(1)}^{2n}.
\]
Note that for $t\geq1$, it follows that $t-(t+1)R_{(1)}>0$ and $t-(t-1)R_{(1)}>0$. Then, we have
\[
\displaystyle
f(-R_{(1)},n)>1-R_{(1)}.
\]Hence, the lemma is proved.
\qed

\medskip

\begin{lem}\label{lem3.2.02}{\normalfont If $\displaystyle\left|\frac{\mathcal{M}_{k-2}(n)}{\mathcal{M}_{k}(n)}\right|_{C_1}<\frac{1}{R_{(1)}}$, then $\displaystyle\left|\frac{\mathcal{M}_{k}(n)}{\mathcal{M}_{k+2}(n)}\right|_{C_1}<\frac{1}{R_{(1)}}$.
}
\end{lem}

\noindent
Proof. From \eqref{eqn-3.2.1}, one has 
\[
\displaystyle
\left|\frac{\mathcal{M}_{k+2}(n)}{\mathcal{M}_{k}(n)}\right|_{C_1}\geq \left|f(q,n)\right|_{C_1}-R_{(1)}^{2n}\left|\frac{\mathcal{M}_{k-2}(n)}{\mathcal{M}_{k}(n)}\right|_{C_1}.
\]
By Lemma \ref{lem3.2.01} and $\displaystyle\left|\frac{\mathcal{M}_{k-2}(n)}{\mathcal{M}_{k}(n)}\right|_{C_1}<\frac{1}{R_{(1)}}$, it follows that
\[
\displaystyle
\left|\frac{\mathcal{M}_{k+2}(n)}{\mathcal{M}_{k}(n)}\right|_{C_1}\geq \left|f(q,n)\right|_{C_1}-R_{(1)}^{2n}\left|\frac{\mathcal{M}_{k-2}(n)}{\mathcal{M}_{k}(n)}\right|_{C_1}>1-R_{(1)}-R_{(1)}^{2n-1}>R_{(1)} \ \text{with} \ n\geq2.
\] On the other hand, the case of $n=1$, has been proved (see \cite[Proposition 3.1]{LMOV}).\qed


\medskip

\begin{lem}\label{lem3.2.03}{\normalfont $\displaystyle f(q, n)$ has no zeros inside $D_1$.
}
\end{lem}

\noindent
Proof. The proof is almost the same as the proof of Theorem \ref{thm-0.2} by the major change being the substitution of 
\begin{equation}\label{eqn-3.2.01}
\displaystyle
\sum^{2n}_{t=0}\mu_tz^t, \ \ 
\displaystyle
\mu_t=\left\{    
\begin{array}{
    cl}
1 &  \text{for}  \   t=0, \\
tr^t  & \text{for}  \ 0<t\leq n, \\
(2n-t)r^t & \text{for}  \  n<t<2n, \\
r^{2n} & \text{for}  \ t=2n \\
\end{array} \right.
\end{equation}
for \eqref{eqn-2.2.1}. We may also prove this lemma by using Lemma \ref{lem3.2.01} and the Rouch\'e theorem. Note that $P_n(q)=f(q,n)+2q^n$ and $P_n(q)$ has no zero inside $D_1$. Then,  we can check that $1-R_{(1)}>2R_{(1)}^n$ with $n\geq2$. Therefore, we obtain 
$
\displaystyle
\left|f(q,n)\right|_{C_1}>2R_{(1)}^n.
$\qed

We are now turning to the proof of Theorem \ref{thm-0.3.1}. We prove that $\mathcal{M}_{k}(n)$ has no zeros inside $D_1$ by induction on $k$.

\noindent

Assume that $\mathcal{M}_{k}(n)$ has no zeros inside $D_1$. By Lemmas \ref{lem3.2.01} and \ref{lem3.2.02}, the following inequality
\[
\displaystyle
\left|f(q,n)\right|_{C_1}\geq R_{(1)}^{2n}\left|\frac{\mathcal{M}_{k-2}(n)}{\mathcal{M}_{k}(n)}\right|_{C_1}
\] 
holds. Since $f(q,n)$ has zeros inside $D_1$, by the Rouch\'e theorem, $\mathcal{M}_{k+2}(n)$ has no zeros inside $D_n$. Hence, by induction on $k$, the proof can be completed.
\qed

\medskip

Since $\widetilde{\mathcal{M}}_{k}(n)$ also satisfies \eqref{eqn-3.1.4}, by the proof of Theorem \ref{thm-0.3.1}, we have the following corollary.

\begin{cor}\label{cor3}{\normalfont
The zeros of the polynomials $\widetilde{\mathcal{M}}_{k}(n)$ and $\mathcal{M}_{k}(n)$ lie in the annulus $\displaystyle \mathcal{A}_1$. 
}
\end{cor}

\medskip


\subsection{Proof of Theorem \ref{thm-0.3}.}
In order to prove Theorem \ref{thm-0.3}, we need the following two lemmas.
\begin{lem}\label{lem3.2.1}{\normalfont On the circles $C_{3}$ and $C_{4}$, one has
\begin{align}
\displaystyle
\left|f(q,3)\right|_{C_3}&>R_{(3)}^4+R_{(3)}^3+R_{(3)}^2,\label{eqn-3.2.2a} \\ 
\left|f(q,4)\right|_{C_4}&>R_{(4)}^6+R_{(4)}^5+R_{(4)}^4+R_{(4)}^3+R_{(4)}^2.\label{eqn-3.2.2b}
\end{align}
}
\end{lem}

\noindent
Proof. By direct calculation, one has
\begin{align}
\displaystyle
R_{(3)}^6-R_{(3)}^5+R_{(3)}^4-4R_{(3)}^3+R_{(3)}^2-R_{(3)}+1&>0,\notag \\
R_{(4)}^8-R_{(4)}^7+R_{(4)}^6-4R_{(4)}^5+3R_{(4)}^4-4R_{(4)}^3+R_{(4)}^2-R_{(4)}+1&>0.\notag
\end{align}
Then 
\[
\begin{split}
\left|f(q,3)\right|_{C_3}&\geq f(-R_{(3)},3)=R_{(3)}^6-R_{(3)}^5+2R_{(3)}^4-3R_{(3)}^3+2R_{(3)}^2-R_{(3)}+1\\
&>R_{(3)}^6-R_{(3)}^5+2R_{(3)}^4-3R_{(3)}-(R_{(3)}^6-R_{(3)}^5+R_{(3)}^4-4R_{(3)}^3+R_{(3)}^2-R_{(3)}+1)\\
&=R_{(3)}^4+R_{(3)}^3+R_{(3)}^2,
\end{split}
\]
\[
\begin{split}
\left|f(q,4)\right|_{C_4}&\geq f(-R_{(4)},4)=R_{(4)}^8-R_{(4)}^7+2R_{(4)}^6-3R_{(4)}^5+4R_{(4)}^4-3R_{(4)}^3+2R_{(4)}^2-R_{(4)}+1\\
&>R_{(4)}^8-R_{(4)}^7+2R_{(4)}^6-3R_{(4)}^5+4R_{(4)}^4-3R_{(4)}^3+2R_{(4)}^2-R_{(4)}+1\\
& \ \ \ \ -(R_{(4)}^8-R_{(4)}^7+R_{(4)}^6-4R_{(4)}^5+3R_{(4)}^4-4R_{(4)}^3+R_{(4)}^2-R_{(4)}+1)\\
&=R_{(4)}^6+R_{(4)}^5+R_{(4)}^4+R_{(4)}^3+R_{(4)}^2. 
\end{split}
\]
\qed

\begin{lem}\label{lem3.2.2}{\normalfont For $n=3,4$, if $\displaystyle\left|\frac{\mathcal{M}_{k-2}(n)}{\mathcal{M}_{k}(n)}\right|_{C_n}<\frac{1}{R_{(n)}}$, then $\displaystyle\left|\frac{\mathcal{M}_{k}(n)}{\mathcal{M}_{k+2}(n)}\right|_{C_n}<\frac{1}{R_{(n)}}$.
}
\end{lem}

\noindent
Proof. 
From \eqref{eqn-3.2.1}, one has 
\[
\displaystyle
\left|\frac{\mathcal{M}_{k+2}(n)}{\mathcal{M}_{k}(n)}\right|_{C_n}\geq \left|f(q,n)\right|_{C_n}-R_{(n)}^{2n}\left|\frac{\mathcal{M}_{k-2}(n)}{\mathcal{M}_{k}(n)}\right|_{C_n}.
\]
By Lemma \ref{lem3.2.1} and $\displaystyle\left|\frac{\mathcal{M}_{k-2}(n)}{\mathcal{M}_{k}(n)}\right|_{C_n}<\frac{1}{R_{(n)}}$, it follows that 
\begin{align}
\displaystyle\left|\frac{\mathcal{M}_{k+2}(3)}{\mathcal{M}_{k}(3)}\right|_{C_3}&>R_{(3)}^4+R_{(3)}^3+R_{(3)}^2-R_{(3)}^5>R_{(3)}\notag \\
\left|\frac{\mathcal{M}_{k+2}(4)}{\mathcal{M}_{k}(4)}\right|_{C_4}&>R_{(4)}^6+R_{(4)}^5+R_{(4)}^4+R_{(4)}^3+R_{(4)}^2-R_{(4)}^7>R_{(4)}.\notag
\end{align}
Hence the lemma is proved.\qed

\medskip

We are now turning to the proof of Theorem \ref{thm-0.3}. By induction on $k$, we prove that $\mathcal{M}_{k}(n)$ has no zeros inside $D_n$ for $n=3,4$.

\noindent

Assume that $\mathcal{M}_{k}(n)$ has no zeros inside $D_n$. Since $R_{(n)}\leq1$, Lemma \ref{lem3.2.1} implies that
\begin{align}
\displaystyle
\left|f(q,3)\right|_{C_3}&>R_{(3)}^4+R_{(3)}^3+R_{(3)}^2>2R_{(3)}^3\geq\left|2q^3\right|_{C_3},\notag \\
\left|f(q,4)\right|_{C_4}&>R_{(4)}^6+R_{(4)}^5+R_{(4)}^4+R_{(4)}^3+R_{(4)}^2>2R_{(4)}^4\geq\left|2q^3\right|_{C_4}.\notag
\end{align}
Since $f(q,n)+2q^n=P_n(q)$, and $P_n(q)$ has no zeros inside $D_n$, we see that $f(q,n)$ has no zeros inside $D_n$, by using the Rouch\'e theorem.

\noindent

On the other hand, by Lemmas \ref{lem3.2.1} and \ref{lem3.2.2}, the following inequality holds.
\[
\displaystyle
\left|f(q,n)\right|_{C_n}\geq R_{(n)}^{2n}\left|\frac{\mathcal{M}_{k-2}(n)}{\mathcal{M}_{k}(n)}\right|_{C_n}
\] 
Again, by the Rouch\'e theorem, $\mathcal{M}_{k+2}(n)$ has no zeros inside $D_n$. By induction on $k$, the proof can be completed.
\qed

\medskip

Since $\widetilde{\mathcal{M}}_{k}(n)$ also satisfies \eqref{eqn-3.1.4}, by the proof of Theorem \ref{thm-0.3}, we have the following corollary which is better than Corollary \ref{cor3} where $n=3,4$.

\begin{cor}\label{cor4}{\normalfont
For $n=3,4$, the zeros of the polynomials $\widetilde{\mathcal{M}}_{k}(n)$ and $\mathcal{M}_{k}(n)$ lie in the annulus $\displaystyle \mathcal{A}_n:=\{q\in\mathbb{C}:R_{(n)}\leq|q|\leq R_{(n)}^{-1}\}$. 
}
\end{cor}

\medskip

\begin{rem}{\normalfont
Note that for the case of $n$ is a specific integer, we can find an inequality to estimate $\left|f(q,n)\right|_{C_n}$ by direct calculation such as Lemma \ref{lem3.2.1}. However, for the general case ($n\geq5$), there are some limitations to this approach to construct an inequality of $\left|f(q,n)\right|_{C_n}$ to check the following conditions.

\begin{enumerate}
\item[(I)] $\left|f(q,n)\right|_{C_n}\geq\left|2q^n\right|_{C_n}$;
\item[(II)] If $\displaystyle\left|\frac{\mathcal{M}_{k-2}(n)}{\mathcal{M}_{k}(n)}\right|_{C_n}<\frac{1}{R_{(n)}}$, then $\displaystyle\left|\frac{\mathcal{M}_{k}(n)}{\mathcal{M}_{k+2}(n)}\right|_{C_n}<\frac{1}{R_{(n)}}$.
\end{enumerate}

}
\end{rem}
\medskip

\section*{Acknowledgments}
I am greatly indebted to Professor Toshiki Matsusaka for many useful discussions, including  the idea of the proof of Theorem \ref{thm-0.2}. I also would like to thank to Professors Takeyoshi Kogiso and Michihisa Wakui who read the paper and made numerous helpful suggestions.

\medskip


\bigskip \bigskip 
\hspace{5cm}
\begin{minipage}[t]{10cm}
Xin Ren
\par 
Department of Mathematics, Faculty of Engineering Science,
Kansai University, Suita-shi, Osaka 564-8680, Japan
\par 
E-mail address: k641241@kansai-u.ac.jp
\par 
\end{minipage}

\end{document}